\documentclass[11pt]{article}
\usepackage[utf8]{inputenc}
\usepackage[
papersize={5.5in, 8.5in},
left=0.35in,
right=0.35in,
top=0.35in,
bottom=0.5in]{geometry}
\usepackage{moresize}
\usepackage{graphicx}
\usepackage{amsthm}
\usepackage{rotating}
\usepackage{amsfonts,amssymb,amsmath}
\usepackage[usenames]{color}
\usepackage[T1]{fontenc}

\usepackage{longtable}
\usepackage{array}
\usepackage{float}
\usepackage{booktabs}
\usepackage{verbatim}
\usepackage[dvipsnames,usenames,table,xcdraw]{xcolor}

\newcommand{\g}{\cellcolor[HTML]{E0E0E0}}

\newcommand\T{\rule{0pt}{2.5ex}}       
\newcolumntype{P}[1]{>{\centering\arraybackslash}p{#1}}

\graphicspath{ {images/} }

\title{\bf 
The chromatic number of $\mathbb{R}^8$ 
is at least 25
}
\author{\bf 
\textcolor[rgb]{0,1.0,0.0}{Jaan Parts} \\
} 
\date{\normalsize \textcolor[rgb]{0,1.0,0.0}{Kazan, Russia, jaan\_parts@.mail.ru}}

\begin{document}

\maketitle

\pagestyle{empty}
\thispagestyle{empty}

\begin{abstract}
The lower bound for the chromatic number of the Euclidean space of dimension 8 is increased to 25.
\end{abstract}

\section{We came in}

A good way to understand something is to try to improve it. Armed with this motto, 
I attempted to tackle the problem of proper coloring of Euclidean spaces in dimensions greater than two.\footnote{
A coloring is called \textit{proper} if any pair of points at (unit) \textit{forbidden distance} from each other takes two different colors. We want to find the minimum possible number of colors, which is called the \textit{chromatic number} $\chi$. There are many variations of this problem, but the simplest and most natural one is to consider it in spaces $\mathbb{R}^n$ of different dimensions $n$. Here, we restrict ourselves to the case $n=8$.}

Before this, I somehow got my bearings in three-dimensional space, but I still cannot imagine four-dimensional space, etc. The beauty of mathematics is that you don't even need to imagine anything to work with it. And once you start working with it, it becomes familiar and not so scary. And this is a good step towards finally being able to imagine and understand it.

Four years ago, I raised the \textit{lower bound} of the chromatic number $\chi(\mathbb{R}^8)$ by one. For a number of reasons (mainly laziness), this result remained unpublished.

The previous record $\chi(\mathbb{R}^8)\ge 19$ was obtained by Matthew Kahle and Birra Taha \cite{kata}. They started with the graph, constructed on 240 vertices of the Gosset polytope $4_{21}$ and having an \textit{independence number} $\alpha = 16$. The idea was to add as many additional vertices to this graph as possible while maintaining the same $\alpha$.

Recently, Aubrey de Grey reminded me of my research in connection with the development of his conjecture (see the Polymath16 project blog \cite{poly}), and I promised to publish it. Upon checking my previous results, I discovered that I could make further progress. As a result, the hunt for graphs yielding record-breaking bounds for $\chi(\mathbb{R}^8)$ dragged on for a month (and it took about the same amount of time to write this article).


\section{Preliminaries}
\paragraph{Approach.}
The lower bound for the chromatic number of a Euclidean space of dimension 8 (and higher) is usually estimated using a finite graph with $v$ vertices and independence number $\alpha$:
\begin{gather*}
\chi\ge \lceil v/\alpha \rceil
\end{gather*}

This relation follows directly from Dirichlet's (pigeonhole) principle: since $\alpha$ is the maximum number of vertices in a graph that can be properly colored with a single color, it is clear that to color the entire graph, one must take at least $\lceil v/\alpha \rceil$ sets of vertices of different colors.

Clearly, this estimate of $\chi$ is quite weak, but we currently lack stronger estimates. All that is required to obtain a lower bound for $\chi(\mathbb{R}^n)$ is a graph for which one can calculate the exact value of the independence number $\alpha$ or at least estimate its upper bound.

An essential element of the approach under consideration is a \textit{solver}, that is, a computer program that calculates the value of $\alpha$ for a graph. On the one hand, the presence of a solver allows us to somewhat improve the estimates of $\chi(\mathbb{R}^n)$. On the other hand, it apparently limits the applicability of this approach to low-dimensional spaces, say $n \le 10$.

The resulting graphs with relatively large $v/\alpha$ are asymmetric, making their manual analysis difficult. To construct such a (good) graph, a promising symmetric core with a sufficiently large number of vertices $v$ and a presumably high $\chi$ is selected, to which, using a good strategy, the vertices are added to maximize the $v/\alpha$ ratio. Apparently, the added vertices do not so much influence the actual $\chi$ value of the graph, as they fill vacant places (by lengthening short independent sets) at a fixed $\alpha$.

Unfortunately, we don't know how to do all these good things. All we can do is try to guess a slightly better solution based on general considerations and experience gained from sifting through a large number of graphs.

\paragraph{Tools.}
To calculate the independence number $\alpha$, we used \textit{exact} solvers that determine the size of the \textit{maximum clique} of the graph complement: \texttt{cliquer} \cite{ost}, \texttt{mcqd} \cite{konc, zhou}, and \texttt{clisat} \cite{san}. To speed up the calculations of $\alpha$ on large graphs, we used the suboptimal solver \texttt{kamis} \cite{lamm}, after which we verified the obtained values using exact solvers.

The upper bound for the chromatic number $\chi$ was determined using the SAT-solver \texttt{yalsat} \cite{bie}.

\paragraph{Notation.}
One of our goals was to provide a simple description of the resulting graphs, so that they demonstrate structure rather than just a jumble of coordinates. We partially succeeded. 

The forbidden distance is set to 4, which excludes fractional values.
The coordinates of individual vertices are specified in full. 
To denote sets of vertices with coordinates that differ in sign and sequence, we use a shorthand notation common to all vertices in the set. We have slightly modified the notation used in the excellent book \cite{con} by John Conway and Neil Sloane to better reflect the peculiarities of our graphs.

In our shorthand notation, the superscript indicates the number of repeating values and the subscript indicates fixed coordinate positions, which are listed without commas. For some values, the superscript may be omitted if it can be unambiguously reconstructed from the context. If the subscript is omitted, all coordinate permutations that do not contradict the fixed positions are used. The signs are specified explicitly. The $\pm$ sign indicates that both signs are used. The prefixes $e$ and $o$ denote an even and odd number of minus signs. The number of vertices in the set is indicated in parentheses.

For example, the coordinates of all 240 vertices of the Gosset polytope $4_{21}$ in the even coordinate system of the $E_8$ lattice are written as follows: $\pm 2^2 \,0^6$ (112), $e1^8$ (128). The additional 49 vertices of 289-vertex Kahle-Taha graph \cite{kata} have the form $\pm 2^4 \,0^4$ (40), $e1^4 \,0^4$ (6), $+3^1 \,o1^7$ (3).

For our purposes, it is convenient to replace the first coordinate $x_1$ with $(2-x_1)$. In this case, the vertices of the Gosset polytope are $0_1 \!\pm\!2^1 \,0^6$ (14), $+2_1\!\pm\!2^2 \,0^5$ (84), $+4_1\!\pm\!2^1 \,0^6$ (14), $+1_1 \,e1^7$ (64), $+3_1 \,o1^7$ (64).

The same notation can be used for individual vertices. For example, $(+2_{1} \,0^7) = (+2_{1}^1 \,0^7) = (2, 0, 0, 0, 0, 0, 0, 0);\;\; (+2_{468} \,0^5) = (+2_{468}^3 \,0^5) = (0, 0, 0, 2, 0, 2, 0, 2)$.

\newpage

\section{Why dont't we do it in the road?}
\paragraph{$\chi(\mathbb{Q}^8)\ge 20$.}
Following Kahle and Taha, we started with the base graph $G_{240}$, constructed on the vertices of the Gosset polytope $4_{21}$.

Without changing $\alpha=16$, we were able to increase the number of additional vertices to 64 (while remaining in the $E_8$ lattice). However, this did not lead to an increase in the lower bound for the chromatic number $\chi(\mathbb{Q}^8)\ge \lceil 304/16 \rceil =19$. The breakthrough was just one vertex away (and this motive haunted us later, driving us to continue the search).

We made progress when we dropped the constraint on $\alpha$ being fixed. With $\alpha=17$, we obtained a 327-vertex graph $G_{327}$ (leaving the $E_8$ lattice). This already yields a bound $\chi\ge 20$.

\paragraph{$\chi(\mathbb{Q}^8)\ge 22$.}
We took a tour of large graphs, believing that a combination of good initial graphs would yield a good final graph (for example, one formed by the Minkowski sum of $G_{240}$ and various cliques).
We investigated a series of graphs formed by successively discarding vertices with lower degree, starting with graphs of several thousand and tens of thousands of vertices. We used the \texttt{kamis} solver to estimate the independence number. At this stage, we obtained graphs with estimates up to $\chi\ge 29$, but \texttt{kamis} is not the exact solver, and we cannot confirm these results. Solvers that produce exact $\alpha$ values perform poorly on graphs with more than 1000 vertices.\footnote{
This threshold is floating and depends on the graph, but it is good enough as a first approximation. By the way, we haven't got a solution for the 1716-vertex graph mentioned in \cite{exis} when estimating the lower bound of $\chi(\mathbb{R}^{12})$, although the \texttt{mcqd} solver was running for over a month.}

As a result, we obtained a 516-vertex graph $G_{516}$ with $\alpha=24$, which gives an estimate $\chi\ge 22$.

\paragraph{$\chi(\mathbb{Q}^8)\ge 24$.}
By filling the vertex sets $\pm2^3\,0^5$ and $\pm 1^8$ (see Table~\ref{tpar}), we obtain a 784-vertex graph $G_{784}$ with $\alpha=34$, which yields $\chi\ge 24$.
If we now remove 64 vertices $+3_1 \,o1^7$, we obtain an intermediate graph $G_{720}$ with $\alpha=33$, which can be supplemented with 48 vertices of the form $\pm 1^4 \,0^4$ without changing $\alpha$, and yielding graph $G_{768}$ also with $\chi\ge 24$.

And these are the last graphs where we were able to operate with relatively large vertex sets. Subsequently, we had to carefully select small vertex sets and individual vertices to add, using the more promising graph $G_{784}$ as a base.

\paragraph{$\chi(\mathbb{Q}^8)\ge 25$.}
One might think that other vertices of the form $\pm 3^1 \!\pm\!1^7$ would be useful for increasing $\chi$. But as it turns out, adding them leads to a rapid increase in $\alpha$, which has the opposite effect. We were only able to make good use of 8 more vertices of this type.
This leads to an interesting phenomenon: asymmetry in only one coordinate. Something similar is observed in sets of the form $\pm 4^1 \!\pm\!2^1 \,0^6$: without changing $\alpha$, only three coordinates are allowed to be occupied by non-zero values.

For a while, we were stuck at 808 vertices. Next jump in $\chi$ was achieved by adding some vertices of the form $\pm 1^4 \,0^4$.

Little by little, we reached the 843-vertex graph $G_{843}$ with $\alpha=34$ and $\chi\ge 25$. It's a shame to stop just a few steps short of another record (851 vertices seemed so close), but it had to be done sometime. 


\section{Graphs}

\paragraph{General information.}
Table~\ref{tpar} lists the vertex subsets and some parameters of the graphs that appear in the text. Subsets containing all possible permutations and sign changes are highlighted in gray.

\begin{table}[h]
{
\caption{Parameters and structure of some graphs.}
\label{tpar}
\small
\smallskip

{
\centering
\begin{tabular}{@{\;}c*{8}{|>{\!\!\centering\arraybackslash}p{8mm}} }
\hline \T
graph			&$G_{327}$ &$G_{347}$ &$G_{516}$ &$G_{720}$ &$G_{768}$ &$G_{784}$ &$G_{818}$ &$G_{843}$ \\
\hline 
\hline \T
vertices ($v$)		&327	&347	&516	&720	&768	&784	&818	&843	\\
edges			&22469	&24459	&45924	&86056	&95296	&94856	&102019	&105180	\\
degree $\ge$	&88 	&91	&148	&229	&193	&155	&160	&67	\\
degree $\le$	&171	&180	&289	&409	&409	&441	&441	&446	\\
\hline \T
$\pm2^1 \;0^7$		&14	&14	&\g 16	&\g 16	&\g 16	&\g 16	&\g 16	&\g 16	\\
$\pm2^3 \;0^5$	&98	&118	&308	&\g 448	&\g 448	&\g 448	&\g 448	&\g 448	\\
$\pm4^1 \pm\!2^1 \;0^6$	&14	&14	&0	&0	&0	&0	&0	&8	\\
$e1^8$			&\g 128	&\g 128	&\g 128	&\g 128	&\g 128	&\g 128	&\g 128	&\g 128	\\
$o1^8$			&0	&0	&0	&\g 128	&\g 128	&\g 128	&\g 128	&\g 128	\\
$+3_1 \;o1^7$		&\g 64	&\g 64	&\g 64	&0	&0	&\g 64	&\g 64	&\g 64	\\
$+3_1 \;e1^7$		&0	&0	&0	&0	&0	&0	&0	&8	\\
$+5_1 \;e1^7$		&0	&0	&0	&0	&0	&0	&0	&8	\\
$e1^4 \;0^4$		&8	&8	&0	&0	&48	&0	&12	&13	\\
$o1^4 \;0^4$		&1	&1	&0	&0	&0	&0	&22	&22	\\
\hline \T
$\alpha$			&17	&18	&24	&33	&33	&34	&34	&34	\\
$v/\alpha$		&19.235	&19.278	&21.5	&21.818	&23.273	&23.059	&24.059	&24.794	\\
$\chi \ge$		&20	&20	&22	&22	&24	&24	&25	&25	\\
$\chi \le$		&21	&22	&26	&24	&24	&27	&27	&27	\\
\hline
\end{tabular}

}
}
\end{table}

\paragraph{Details.}

The graphs $G_{327}$ and $G_{347}$ with chromatic number $\chi \ge 20$ include 240 vertices of the Gosset polytope (we use the form with biased first coordinate) and all 64 vertices $-1_1 \,e1^7$. The remaining vertices are of the form $0_1 \!\pm\! 2^3 \,0^4$ and $+1_1 \!\pm\! 1^3\, 0^4$. Here is a list of coordinates of the additional 23 vertices of the graph $G_{327}$:

\vspace{-2mm}
\begin{table}[H]
{
\small
\begin{tabular}{*{2}{p{56mm}}}
\centering
\begin{tabular}{*{8}{>{\!\!\!\!\raggedleft\arraybackslash}p{2mm}}}
0	&$-$2	&$-$2	&2  	&0  	&0  	&0  	&0  \\
0	&$-$2	&$-$2	&$-$2	&0  	&0  	&0  	&0  \\
0	&$-$2	&0  	&0  	&$-$2	&2  	&0  	&0  \\
0	&$-$2	&0   	&0  	&$-$2	&$-$2	&0  	&0  \\
0	&$-$2	&0  	&0  	&0  	&0  	&$-$2	&2  \\
0	&$-$2	&0  	&0  	&0  	&0  	&$-$2	&$-$2 \\
0	&0	    &$-$2	&0  	&$-$2	&0  	&2  	&0  \\
0	&0	    &$-$2	&0  	&$-$2	&0  	&$-$2	&0  \\
0	&0	    &$-$2	&0  	&0  	&$-$2	&0  	&2  \\
0	&0	    &$-$2	&0  	&0  	&$-$2	&0  	&$-$2 \\
0	&0  	&0  	&$-$2	&$-$2	&0  	&0  	&2  \\
0	&0	    &0  	&$-$2	&$-$2	&0  	&0  	&$-$2 \\
\end{tabular}
&
\centering
\begin{tabular}{*{8}{>{\!\!\!\!\raggedleft\arraybackslash}p{2mm}}}
0	&0  	&0  	&$-$2	&0  	&$-$2	&2  	&0  \\
0	&0	    &0  	&$-$2	&0  	&$-$2	&$-$2	&0  \\
1	&0  	&$-$1	&$-$1	&1  	&0  	&0  	&0  \\
1	&0  	&$-$1	&$-$1	&0  	&0  	&1  	&0  \\
1	&0  	&$-$1	&1  	&$-$1	&0  	&0  	&0  \\
1	&0  	&$-$1	&1  	&0  	&$-$1	&0  	&0  \\
1	&0  	&$-$1	&0  	&$-$1	&0  	&0  	&1  \\
1	&0  	&$-$1	&0  	&0  	&1  	&0  	&$-$1 \\
1	&0  	&1  	&0  	&$-$1	&$-$1	&0  	&0  \\
1	&0  	&1  	&0  	&1  	&0  	&$-$1	&0  \\
1	&0  	&0  	&$-$1	&$-$1	&1  	&0  	&0  \\
\\
\end{tabular}
\end{tabular}
}
\end{table}
\vspace{-2mm}

Graph $G_{347}$ is formed from $G_{327}$ by adding 20 vertices $0_1 \pm\!2^3 \,0^4$, which we will not show here, since $G_{347}$ is not of particular interest.

Graph $G_{516}$ includes all vertices of the form $\pm 2^1 \,0^7$ (16), $+2_1 \pm\!2^2 \,0^5$ (84), $e1^8$ (128), $+3_1 \,o1^7$ (64) and vertices of the form $0_1 \!\pm\!2^3 \,0^4$ (224), except for following 56 vertices: $\pm2_{234} \,0^5$ (8), $\pm2_{256} \,0^5$ (8), $\pm2_{278} \,0^5$ (8), $\pm2_{357} \,0^5$ (8), $\pm2_{367} \,0^5$ (8), $\pm2_{457} \,0^5$ (8), $\pm2_{468} \,0^5$ (8).

The graphs $G_{720}$ and $G_{784}$ are completely specified by the subsets shown in Table~\ref{tpar}.

Graph $G_{720}$ is a base subgraph for the graphs $G_{768}$, $G_{784}$, $G_{818}$, and $G_{843}$ described below. It has symmetry order 10321920 and 3 orbits containing all vertices of the form $\pm 2^1\,0^7$ (16), $\pm 2^3\,0^5$ (448), and $\pm 1^8$ (256).

Graph $G_{768}$ includes additional vertices $e1^4 \,0^4$ (48), namely, all vertices $e1_{1234} \,0^4$ (8), $e1_{1256} \,0^4$ (8), $e1_{1278} \,0^4$ (8), $e1_{3456} \,0^4$ (8), $e1_{3478} \,0^4$ (8), $e1_{5678} \,0^4$ (8). 

Graph $G_{784}$ includes additional vertices $+3_1 \,o1^7$ (64).


Graph $G_{818}$ is obtained from $G_{784}$ by adding 34 vertices of the form $\pm 1^4 \,0^4$, including three quadruples of vertices $+1_1 \,e1^3 \,0^4$ (12), namely, $+1_1 \,e1_{368} \,0^4$ (4), $+1_1 \,e1_{458} \,0^4$ (4), $+1_1 \,e1_{467} \,0^4$ (4), and 11 pairs of vertices of the form $\pm (o1^4 \,0^4)$ (half of the vertices are shown below, the rest are obtained by changing the sign of all coordinates):

\vspace{-3mm}
\begin{table}[H]
{
\small
\begin{tabular}{*{2}{p{56mm}}}
\centering
\begin{tabular}{*{8}{>{\!\!\!\!\raggedleft\arraybackslash}p{2mm}}}
0  	&0  	&0  	&0  	&1  	&1  	&1  	&$-$1  	\\
0  	&0  	&0  	&0  	&1  	&1  	&$-$1  	&1  	\\
0  	&0  	&1  	&$-$1  	&0  	&0  	&1  	&1  	\\
0  	&0  	&$-$1  	&1  	&0  	&0  	&1  	&1  	\\
0  	&0  	&1  	&1  	&1  	&$-$1  	&0  	&0  	\\
0  	&0  	&1  	&1  	&$-$1  	&1  	&0  	&0  	\\
\end{tabular}
&
\centering
\begin{tabular}{*{8}{>{\!\!\!\!\raggedleft\arraybackslash}p{2mm}}}
0  	&1  	&0  	&1  	&$-$1  	&0  	&1  	&0  	\\
0  	&$-$1  	&0  	&1  	&1  	&0  	&1  	&0  	\\
1  	&0  	&$-$1  	&0  	&0  	&1  	&0  	&1  	\\
$-$1  	&0  	&1  	&0  	&0  	&1  	&0  	&1  	\\
0  	&1  	&$-$1  	&0  	&0  	&1  	&0  	&1  	\\
\\
\end{tabular}
\end{tabular}
}
\end{table}
\vspace{-4mm}

This is sufficient to reach a lower bound $\chi(\mathbb{R}^8)\ge 25$, but we provide a full list of vertices of the graph $G_{843}$ with the maximum $v/\alpha$ ratio we achieved, in case anyone wants to go further.

Graph $G_{843}$ is obtained by adding to graph $G_{818}$ the vertices of the form $\pm 4^1 \!\pm\!2^1 \,0^6$ (8) with nonzero values of the first three coordinates: 

\vspace{-3mm}
\begin{table}[H]
{
\small
\begin{tabular}{*{2}{p{56mm}}}
\centering
\begin{tabular}{*{8}{>{\!\!\!\!\raggedleft\arraybackslash}p{2mm}}}
4  	&2  	&0  	&0  	&0  	&0  	&0  	&0     	\\
4  	&$-$2  	&0  	&0  	&0  	&0  	&0  	&0  	\\
4  	&0  	&2  	&0  	&0  	&0  	&0  	&0  	\\
4  	&0  	&$-$2  	&0  	&0  	&0  	&0  	&0  	\\
\end{tabular}
&
\centering
\begin{tabular}{*{8}{>{\!\!\!\!\raggedleft\arraybackslash}p{2mm}}}
2  	&4  	&0  	&0  	&0  	&0  	&0  	&0  	\\
$-$2  	&4  	&0  	&0  	&0  	&0  	&0  	&0  	\\
2  	&0  	&4  	&0  	&0  	&0  	&0  	&0  	\\
2  	&0  	&$-$4  	&0  	&0  	&0  	&0  	&0  	\\
\end{tabular}
\end{tabular}
}
\end{table}
\vspace{-4mm}

\noindent
vertex sets $+3_1 \,e1^7$ (8), $+5_1 \,e1^7$ (8), differing only in the first coordinate (only the first set is shown below): 

\vspace{-3mm}
\begin{table}[H]
{
\small
\begin{tabular}{*{2}{p{56mm}}}
\centering
\begin{tabular}{*{8}{>{\!\!\!\!\raggedleft\arraybackslash}p{2mm}}}
3  	&1  	&1  	&1  	&1  	&1  	&1  	&1  	\\
3  	&1  	&1  	&1  	&$-$1  	&$-$1  	&$-$1  	&$-$1  	\\
3  	&1  	&$-$1  	&$-$1  	&$-$1  	&$-$1  	&1  	&1  	\\
3  	&1  	&$-$1  	&$-$1  	&1  	&1  	&$-$1  	&$-$1  	\\
\end{tabular}
&
\centering
\begin{tabular}{*{8}{>{\!\!\!\!\raggedleft\arraybackslash}p{2mm}}}
3  	&$-$1  	&1  	&$-$1  	&1  	&$-$1  	&1  	&$-$1  	\\
3  	&$-$1  	&1  	&$-$1  	&$-$1  	&1  	&$-$1  	&1  	\\
3  	&$-$1  	&$-$1  	&1  	&$-$1  	&1  	&1  	&$-$1  	\\
3  	&$-$1  	&$-$1  	&1  	&1  	&$-$1  	&$-$1  	&1  	\\
\end{tabular}
\end{tabular}
}
\end{table}
\vspace{-4mm}

\noindent
and one additional vertex:

\vspace{-3mm}
\begin{table}[H]
{
\small
\begin{tabular}{*{2}{p{56mm}}}
\centering
\begin{tabular}{*{8}{>{\!\!\!\!\raggedleft\arraybackslash}p{2mm}}}
0  	&1  	&0  	&$-$1  	&$-$1  	&0  	&1  	&0  	\\
\end{tabular}
&
\centering
\begin{tabular}{*{8}{>{\!\!\!\!\raggedleft\arraybackslash}p{2mm}}}
\\
\end{tabular}
\end{tabular}
}
\end{table}
\vspace{-5mm}

\paragraph{Computation time.}
The efficiency of programs for calculating $\alpha$ can be compared using the graphs considered above. For example, the solvers \texttt{cliquer}, \texttt{mcqd}, and \texttt{clisat} complete the graph $G_{516}$ in 2650, 94, and 64 seconds, respectively. By comparison, the function \texttt{FindIndependentSet} in \texttt{Wolfram Mathematica 13.3} takes about 11 hours. 

The solvers \texttt{mcqd} and \texttt{clisat} show comparable results on asymmetric graphs: both solvers take about 200 minutes to complete the graph $G_{843}$. On highly symmetric graphs, \texttt{clisat} performs better. We ran \texttt{clisat} with the following parameters (which seem to work faster): vertex ordering $=1$ (DEG-SORT), heuristic $=1$ (AMTS heuristic).
The order of vertices when running \texttt{mcqd} and \texttt{clisat} has little effect on the computation time, since solvers repeatedly apply the built-in suboptimal sorting procedure.

The calculations were performed on a computer with a processor AMD Ryzen 9 7940HS, 4 GHz, 8 cores/16 threads, and 64 GB of RAM, in a single thread under full CPU load.

\newpage
\section{Isn’t this where…}

Initially, we believed that by using only the vertices of the $E_8$ lattice, we could obtain lower bounds for $\chi$ that were close to optimal for a limited number of vertices. If this is true, then the search area is significantly narrowed. This hypothesis was partially confirmed. However, vertices that extend beyond the lattice also make a significant contribution to the estimate of $\chi$ in the graphs we considered. Therefore, we are inclined to believe that such vertices are essential for record-breaking constructions.

There are other questions to which we have only partial answers:

$\bullet$  How can we guess a good base subgraph?

$\bullet$ What is the optimal vertex addition strategy to obtain a relatively high lower bound for $\chi$?

$\bullet$ How far can we push the lower bound for $\chi$ using known approaches?

$\bullet$ How symmetrical should a good resulting graph be?

In our experience, a good base subgraph has good symmetry. It then grows with subgraphs of progressively less symmetry and order, until we reach scattered vertices that are taken without any system.

But perhaps the asymmetry of the resulting graph indicates that we made an unfortunate choice when adding the next batch of vertices, thereby blocking the possibility of adding more vertices.

One can spend a long time occupied with the fascinating task of finding better bounds for the chromatic number of 
$\mathbb{R}^8$. 
Clearly, the potential for further improvement is far from exhausted.
And there are also other dimensions and other estimation methods (see, for example, \cite{cher, exis}).
In fact, we have already obtained the estimate $\chi(\mathbb{Q}^8)>25$ using a slightly different approach, but that is another story.

\end{document}